\documentclass[a4paper,11pt]{article}
\usepackage{amssymb,rotating,graphics,epsfig,float,graphics}
\usepackage{amsmath,verbatim,a4wide}
\usepackage[english]{babel}
\usepackage{amsmath,amsfonts,amssymb,pifont}
\usepackage{textcomp}
\usepackage{latexsym}
\usepackage{theorem,srcltx}
\usepackage{graphicx,float}
\usepackage{mathpazo} 
\usepackage{longtable}
\usepackage{txfonts}

\usepackage{theorem}
\theoremstyle{break} 
{\theorembodyfont{\rmfamily} }
{\theorembodyfont{\rmfamily} }
{\theorembodyfont{\rmfamily} }
{\theorembodyfont{\rmfamily} }

\begin{document}

\author{Ferdinand Verhulst\\
University of Utrecht, Department of mathematics \\
PO Box 80.010, 3508 TA Utrecht, The Netherlands } 

\title{Evolution to symmetry}
\date{ Accepted for publication in {\em Symmetry} 2021}
\maketitle
\begin{abstract} 
A natural example of evolution can be described by a time-dependent two 
degrees-of-freedom Hamiltonian. 
We choose the case where initially the Hamiltonian derives from a general cubic potential, the linearised system 
has frequencies 1 and $\omega >0$.  The time-dependence  produces slow evolution to discrete (mirror) symmetry 
in one of the degrees-of-freedom. This changes the dynamics drastically depending on the frequency ratio $\omega$ 
and the timescale of evolution. We analyse the cases $\omega = 1, 2, 3$ where the ratio's 1,2 turn out to be the most 
interesting. In an initial phase we find 2 adiabatic invariants with changes near the end of evolution.  A remarkable 
feature is the vanishing and emergence of normal modes, stability changes and strong changes of the velocity 
distribution in phase-space.  The problem  is inspired by the dynamics of axisymmetric, rotating galaxies 
that evolve slowly to mirror symmetry with respect to the galactic plane, the model formulation is quite general. 
\end{abstract} 

MSC classes:	37J20, 37J40, 34C20, 58K70, 37G05, 70H33, 70K30, 70K45\\


\section{Introduction} 
There are many dynamical systems with evolutionary aspects (time-dependence) but the necessary theory for 
understanding them is still restricted. The main reason for 
this is that understanding steady state (autonomous) systems  is already a formidable task. 
Physical examples are the evolution to spherical structures in nature and on solar system scale planetary 
satellite systems evolving by long-term tidal forces to more symmetric orbits (for references see \cite{VH08}). \\ 
In \cite{HV97} a cartoon problem is considered of the form: 
\begin{equation}  \label{cart1} 
\ddot{x} + x = a(\varepsilon t)x^2, 
\end{equation} 
where a dot means differentiation with respect to time, $\varepsilon$ is a positive small pameter. The function 
$a(\varepsilon t)$  is monotonically decreasing from the value 1 to zero.  Eq. \eqref{cart1} models slow evolution 
to symmetric dynamics; although the equation is simple it shows already a relation with dissipative systems and the 
construction of a global adiabatic invariant. \\ 
In \cite{VH08}  a 2 degrees-of-freedom (dof) system is considered with a cubic potential that is discrete symmetric in one 
of the positions ($q_2$) and asymmetric in the other position ($q_1$). The Hamiltonian is: 
\begin{equation} \label{cart2} 
H= \frac{1}{2}(p_1^2+4q_1^2) + \frac{1}{2}(p_2^2 +q_2^2) + a(\varepsilon t) (2q_1^3 +q_1q_2^2). 
\end{equation}
As before the asymmetric term vanishes slowly, the system evolves to a 
2 -dimensional harmonic oscillator. Again a relation with dissipative systems can be established, 2 adiabatic invariants 
can be found. The system displays overall dynamics that keeps some information from its asymmetric past. 

 \subsection{The collisionless Boltzmann equation} 
 The (to a good approximation) lack of collisions in stellar systems raises questions on the statistical mechanics of 
 these systems. The collisionless Boltzmann  or Liouville equation describes the distribution of $n$ particles, stars in this 
 case, in $6n$ dimensional phase-space. A good description with examples can be found in \cite{BT} ch.~4. The equation 
 for the distribution of particles $f(t, {\bf x}, {\bf v})$ where ${\bf x}$ indicates the position, ${\bf v}$ the velocity, is the 
 continuity equation in  $6n$ dimensional phase-space with clearly the assumption that no particle escapes, is destroyed 
 or is created in the system. With Lagrangian derivative $d/dt$ the Liouville equation is $df/dt =0$. If the collective 
 gravitational potential ruling the dynamics of the system is $\Phi(t, {\bf x})$, the Liouville equation becomes explicitly: 
 \begin{equation} \label{L} 
 \frac{\partial f}{\partial t} + \sum_{i=1}^3\left(v_i  \frac{\partial f}{\partial x_i} - \frac{\partial \Phi}{\partial x_i}
 \frac{\partial f}{\partial v_i}\right)=0. 
 \end{equation}  
 The characteristics of this first order partial differential equation are given by the Hamiltonian equations of motion. 
 The solutions of the Hamiltonian system produce according to Monge the solutions of the Liouville equation as they 
 represent the geometric sets where the solutions of the Liouville equation are constant. This procedure is very effecrive 
 if we find not only solutions but integrals of motion that contain sets of solutions. Assuming 
 for instance a time-independent potential and axi-symmetry, we have already 2 independent integrals of motion, energy 
 $E$ and angular momentum $L$ with respect to the axis of rotation. Any differential function $f$ of $E$ and $L$ will 
 satisfy the Liouville equation. It was noted very early that such solutions will produce velocity distributions that 
 are symmetric perpendicular to and in the direction of the rotation axis. This does not agree with observations in our 
 own galaxy. This triggered off a long search for ``a third integral of the galaxy'' that ended in the 1960's when it was shown 
 in Arnold-Moser theory that such 3rd integrals in general do not exist.\\
 One of the aims of our study is to investigate whether evolution from an asymmetric state to an integrable symmetric 
 dynamical state 
 might influence the velocity distributions as if the system ``remembers'' its asymmetric past.\\
 From \cite{BT} ch.~4 we consider the dynamics in cylindrical coordinates $R, \phi, z$ with $\dot{R}=v_R, 
 \dot{\phi}= v_{\phi}/R, \dot{z}=v_z$; we have with potential $\Phi(R, \phi, z)$ the equations of motion: 
 \begin{equation} \label{eqscyl} 
 \dot{v}_R = - \frac{\partial \Phi}{\partial R} + \frac{v_{\phi}^2}{R}, \dot{v}_{\phi}= - \frac{1}{R}\frac{\partial \Phi}{\partial \phi} 
 - \frac{v_Rv_{\phi}}{R}, \dot{v}_z= - \frac{\partial \Phi}{\partial z}. 
 \end{equation}

 \subsection{Axi-symmetric models} 
 In the sequel we will consider models that are inspired by axisymmetric rotating galaxies. 
 One can think of  disk galaxies or rotating flattened elliptical galaxies.  The axisymmetry expessed by 
 \begin{equation} \label{axi}
 \frac{\partial \Phi}{\partial \phi}  =0,
 \end{equation}
 produces the 
angular momentum integral $J$ (in \cite{BT} called $L_z$) enabling us to reduce 3-dimensional motion to 2 dimensions. 
We have with eqs.~\eqref{eqscyl},  \eqref{axi} 
\[ \dot{v}_{\phi}=  - \frac{v_Rv_{\phi}}{R}\,\, {\rm or} \,\, R \ddot{\phi} + 2 \dot{R} \dot{\phi}=0, \]
leading to the angular momentum integral 
\begin{equation} \label{mom}
R^2 \dot{\phi} = J. 
\end{equation}
The equations of motion \eqref{eqscyl} can be written as 
\[ \ddot{R} =  \frac{J^2}{R^3} - \frac{\partial \Phi}{\partial R} . \]
In the equatorial plane that 
is perpendicular to the axis of rotation we have circular orbits at $R=R_0$ where the rotation speed matches constant 
angular momentum: 
\[   \frac{J^2}{R_0^3} = \frac{\partial \Phi}{\partial R}(R_0,0) . \] 
We will expand the Hamiitonian around the circular orbits where $q_1$ corresponds with the radial direction, $q_2$ with 
expansion in the $z$-direction. 
In most models, the potential $\Phi$ is assumed to be symmetric with respect to the equatorial plane. In cylindrical 
coordinates $R, z$ with $z$ in the direction of the rotation axis and $z=0$ corresponding with the equatorial plane we will put 
in a final stage of evolution $\Phi = \Phi(R, z^2)$. The evolution towards this symmetric state is caused by mechanisms 
unknown to us, maybe contraction combined with rotation or dynamical friction plays a part. We propose to avoid the 
speculative description of complicated 
mechanisms by introducing a function of time slowly destroying the asymmetric potential. \\
Around the circular orbits in the galactic plane
 we find epyciclic orbits in the $R$-direction nonlinearly coupled to bounded  vertical motion in the $z$-direction. An early 
 study of such orbits in a steady state galaxy model is \cite{O62}, for a systematic evaluation of the theory see \cite{BT}. 
 A detailed analysis of the orbits can be found in \cite{FV79}.

\section{A two degrees-of-freedom model with evolution} 
Consider the time-dependent two dof Hamiltonian: 
\begin{equation} \label{Ham} 
H = \frac{1}{2}(\dot{q}_1^2 +q_1^2) + \frac{1}{2}(\dot{q}_2^2 +\omega^2 q_2^2) -(\frac{1}{3}a_1q_1^3 +a_2q_1q_2^2) 
- \alpha(\delta t) (\frac{1}{3}a_3q_2^3 + a_4q_1^2q_2). 
\end{equation} 
The epicyclic frequency has beem scaled to 1, the vertical frequency is $\omega$. 
The function $ \alpha(\delta t)$ is continuous and monotonically decreasing from $\alpha(0)=1$ to zero; in examples 
we take $\alpha(\delta t)=  e^{- \delta t} $. 
If $a_1=1, a_2=-1, a_3=a_4=0$ we have the famous H\'enon-Heiles problem, \cite{HH64}.\\ 
The quadratic part is called $H_2$, the cubic part $H_3$. We assume $0 < \delta \ll 1$, the coefficients 
$a_1, \ldots, a_4$ are free parameters with $a_4 \neq 0$, $\omega$ represents the frequency ratio of the two dof 
with prominent resonances $\omega =1, 2, 3$.   The values $\omega$ can take depend on the galactic potential 
constructed. An example describing an axisymmetric rotating oblate galaxy can be found in \cite{BT} eq. (3-50), leading 
to frequency ratios given by eq. (3-56). 

The terms with coefficients $a_1, a_2$ are discrete symmetric in $q_2$, 
de terms with $a_3, a_4$ are not symmetric in $q_2$ but the asymmetry vanishes as $t \rightarrow \infty$.  
So in a model of a rotating galaxy $q_2$ corresponds with $z$. 

The dynamical system induced by \eqref{Ham} is not reversible as time-independent Hamiltonians are, 
the main question of interest is then whether after a long time the system induced by the Hamiltonian shows traces 
of the original asymmetry.\\
If $\delta =0$ the origin of phase-space is Lyapunov-stable, the energy manifold is bounded in an $O(\varepsilon)$ 
neighbourhood of the origin. 
To make the local analysis more transparent we rescale the coordinates $q_1 = \varepsilon \bar{q}_1$ etc. Dividing by 
$\varepsilon^2$ and leaving out the bars we obtain the equations of motion: 
\begin{eqnarray} \label{eqs}
\begin{cases}
\ddot{q}_1 + q_1 & =   \varepsilon (a_1q_1^2 +a_2q_2^2) + \varepsilon  \alpha(\delta t) 2a_4q_1q_2,\\ 
\ddot{q}_2 + \omega^2 q_2 & =   \varepsilon 2a_2q_1q_2 + \varepsilon  \alpha(\delta t) (a_3q_2^2 +a_4q_1^2).
\end{cases} 
\end{eqnarray} 
Because of the localisation near the origin of phase-space we assume that $\varepsilon$ is a small positive parameter. The 
parameter $\delta$ will be very small, for instance if we want to study the influence of dynamical friction in a galaxy. 
Choosing $\delta = \varepsilon^n$, this implies that for $n=1$ and so 
$ \varepsilon = \delta$ we will consider a very small neighbourhood of the origin. On choosing $ \varepsilon = \surd \delta$
or $n=2$, the neighbourhood will be larger. Again a larger neigbourhood is obtained for $n=3$, these different cases 
complicate the analysis and will produce different local dynamics. 

\section{First order averaging-normalisation} \label{sec3} 
To characterise the dynamics induced by Hamiltonian \eqref{Ham} we will use averaging-normalisation, see \cite{SVM} 
or an introduction in \cite{FV85}. We transform to slowly varying polar coordinates $r, \psi$ by:
\begin{equation} \label{trans} 
q_1=r_1 \cos(t+ \psi_1), \dot{q}_1= -r_1 \sin(t+ \psi_1), q_2=r_2 \cos( \omega t+ \psi_2), 
\dot{q}_2= - \omega r_2 \sin(\omega t+ \psi_2), 
\end{equation}  
leading to the slowly varying system: 
\begin{eqnarray} \label{slow}
\begin{cases} 
\dot{r}_1 & = - \varepsilon \sin(t+ \psi_1) \left(a_1 r_1^2 \cos^2(t+ \psi_1) + a_2 r_2^2 \cos^2( \omega t + \psi_2)\right) - \\ &
 \hspace{5mm}
\varepsilon  \alpha(\delta t) \sin(t+ \psi_1) 2a_4 r_1 \cos(t+ \psi_1) r_2 \cos(\omega t + \psi_2),\\ 
\dot{\psi}_1 & = - \varepsilon \frac{\cos(t+ \psi_1)}{r_1} \left(a_1 r_1^2 \cos^2(t+ \psi_1) + a_2 r_2^2 \cos^2( \omega t +
 \psi_2)\right) - \\ &  \hspace{5mm} \varepsilon  \alpha(\delta t) \cos(t+ \psi_1) 2a_4 \cos(t+ \psi_1) r_2 \cos(\omega t + 
 \psi_2),\\

\dot{r}_2 & = - \frac{\varepsilon}{\omega} \sin(\omega t +\psi_2) 2a_2 r_1 \cos(t+ \psi_1)r_2 \cos(\omega t + \psi_2) - \\ &  
\hspace{5mm} 
\frac{\varepsilon}{\omega}  \alpha(\delta t)  \sin(\omega t +\psi_2) \left(a_3r_2^2 \cos^2(\omega t + \psi_2) +a_4r_1^2 \cos^2 
(t+ \psi_1)\right),\\ 

\dot{\psi}_2 & = - \varepsilon \frac{\cos(\omega t+ \psi_2)}{\omega r_2} 2a_2 r_1 \cos(t+ \psi_1)r_2 \cos(\omega t + \psi_2) - 
\\ &  \hspace{5mm}  - \varepsilon \frac{\cos(\omega t+ \psi_2)}{\omega r_2}  \alpha(\delta t) \sin(\omega t +\psi_2)
 \left(a_3r_2^2 \cos^2(\omega t + \psi_2) +a_4r_1^2 \cos^2 (t+ \psi_1)\right). 
\end{cases} 
\end{eqnarray} 
Near the normal modes $r_1=0$ and $r_2=0$ we have to use a different coordinate transformation. 

We put $\tau = \delta t$ and treat $\tau$ as a new variable. It will also be useful to introduce the actions $E_1, E_2$ by: 
\begin{equation} \label{actions}
E_1= \frac{1}{2}(\dot{q}_1^2 + q_1^2)= \frac{1}{2}r_1^2,\,E_2= \frac{1}{2}(\dot{q}_2^2 + \omega^2 q_2^2)= 
\frac{\omega^2}{2}r_2^2. 
\end{equation}

 As we shall see in subsequent sections, for each choice of $\omega \geq 1$ the average of the terms in system \eqref{slow} 
with coefficients $a_1, a_2$ vanish, so we can use the near-identity transformation \eqref{nearid} of section \ref{app}. 
The implication is {\em that to first order in} $\varepsilon$ and on time intervals of size $1/ \varepsilon$ system \eqref{eqs} 
is described by the intermediate normal form equations:
\begin{eqnarray} \label{eqs2}
\begin{cases}
\ddot{q}_1 + q_1 & =    \varepsilon  \alpha(\delta t) 2a_4q_1q_2,\\ 
\ddot{q}_2 + \omega^2 q_2 & =    \varepsilon  \alpha(\delta t) (a_3q_2^2 +a_4q_1^2).
\end{cases} 
\end{eqnarray} 
We recognize the presence of the $q_2, \dot{q}_2$ normal mode solution as $q_1= \dot{q}_1=0$ satisfies the system; 
this can be checked by using a coordinate system different from polar coordinates.
. 
As the nonlinear terms are homogeneous in the coordinates we can remove the time-dependent term by a 
transformation involving $ \alpha(\delta t)$. For instance if $ \alpha(\delta t) = e^{- \delta t}$ we put  
$q_1 =  \alpha(\delta t) z_1, q_2 = e^{\delta t} z_2$ (note that such a transformation exists for any positive sufficiently 
differentiable function of time that decreases monotonically to zero). System \eqref{eqs2} transforms to: 
\begin{eqnarray} \label{eqs3}
\begin{cases}
\ddot{z}_1 + z_1 & =    -2 \delta \dot{z}_1 - \delta^2 z_1 + \varepsilon  2a_4z_1z_2,\\ 
\ddot{z}_2 + \omega^2 z_2 & =     -2 \delta \dot{z}_2 - \delta^2 z_2 + \varepsilon  (a_3z_2^2 +a_4z_1^2).
\end{cases} 
\end{eqnarray} 
So the time-dependence removing the asymmetry transforms to a dissipative system with friction coeficient $2 \delta$. \\

We will average the righthand sides of system 
\eqref{slow} over $t$ keeping $r_1, r_2, \psi_1, \psi_2, \tau$ fixed. We have $\dot{\tau}= \delta$, so to match the 
size of the other equations of system \eqref{slow} we choose $\delta = \varepsilon^n$ with a suitable choice of 
$n \geq 1$; for simplicity we restrict $n$ to natural numbers. 
As stated above, by first-order averaging the 
terms with coefficients $a_1, a_2$ vanish, we can use system \eqref{eqs2}. 
The subsequent averaging results depend strongly on the choice of $\omega$.  To 
make the calculations more explicit we put in the sequel 
\begin{equation} \label{alfa} 
\alpha(\delta t) = e^{- \delta t},\, \delta = \varepsilon^n (n= 1, 2, \cdots). 
\end{equation} 
On choosing polynomial decrease 
of $\alpha(\delta t)$ we would have slower decrease with as a consequence that we have to retain more small perturbation 
terms. 

\section{The $1:2$ resonance} \label{sec4}
The prominent case for 2 dof systems is the $1:2$ resonance ($\omega =2$). We analyse the system for different 
choices of $\delta$. 
\subsection{First order averaging} 
Averaging system  \eqref{eqs2} we find: 
\begin{eqnarray} \label{ave1}
\begin{cases} 
\dot{r}_1 & = - \varepsilon e^{- \tau} \frac{a_4}{2} r_1r_2 \sin \chi, \,\dot{\psi}_1 = - \varepsilon e^{- \tau} \frac{a_4}{2} 
r_2 \cos \chi,\\ 
\dot{r}_2 & = \varepsilon  e^{- \tau} \frac{a_4}{8} r_1^2 \sin \chi,\,\dot{\psi}_2 =  - \varepsilon e^{- \tau} 
\frac{a_4}{8} \frac{r_1^2}{r_2} \cos \chi. 
\end{cases} 
\end{eqnarray} 
with combination angle $\chi = 2 \psi_1 - \psi_2$ and for $\chi$ the equation: 
\begin{equation} \label{comb12} 
\frac{d \chi}{dt} = \varepsilon a_4 e^{- \tau} (-r_2+ \frac{1}{8} \frac{r_1^2}{r_2}) \cos \chi. 
\end{equation} 
Remarkably, system \eqref{ave1} admits families of solutions with constant amplitude on intervals $O(1/ \varepsilon)$ if: 
\begin{equation} \label{resman1}
 \chi = 0, \pi, r_1^2 = 8r_2^2. 
 \end{equation}
 If $\delta = \varepsilon$, the corresponding phases are slowly decreasing at the rate $\exp{(- \tau)}$. 

If we choose $\delta = \varepsilon^n, n \geq 2$,, the amplitudes $r_{1, 2}$ and phases $\psi_{1, 2}$ will be constant 
with error $O(\varepsilon)$ on intervals of time  of size $1/ \varepsilon$.

System \eqref{ave1} admits 2 time-independent  integrals of motion: 
\begin{equation} \label{intave1}
\frac{1}{2}r_1^2 + 2r_2^2= E_0 \,\,{\rm {and}}\,\, a_4r_1^2r_2 \cos \chi = I_3,  
\end{equation} 
with constants $E_0, I_3$; In the original coordinates we have: 
\[ \frac{1}{2}(\dot{q}_1^2 +q_1^2) + \frac{1}{2}(\dot{q}_2^2 + 4q_2^2)=E_0,\,\,  
a_4(q_1^2q_2- \dot{q}_1^2q_2 + 2q_1\dot{q}_1\dot{q}_2)=I_3. \] 
The solutions and integrals (adiabatic invariants) of system \eqref{ave1} have the error estimate $O(\varepsilon)$ on 
time intervals of size $1/ \varepsilon$. 
On this long interval of time and longer ones we expect the terms $O(\varepsilon^2)$ to play a part as the solutions of 
system \eqref{ave1} 
with coefficients $a_3, a_4$ will vanish and other terms of system \eqref{eqs} will become important. 

\subsection{Second order averaging} 
Second order averaging of system \eqref{slow} produces the system: 
\begin{eqnarray} \label{ave2}
\begin{cases} 
\dot{r}_1 & = - \varepsilon e^{- \tau} \frac{a_4}{2} r_1r_2 \sin \chi, \,
\dot{r}_2  = \varepsilon  e^{- \tau} \frac{a_4}{8} r_1^2 \sin \chi,\\ 

\dot{\psi}_1 & = - \varepsilon e^{- \tau} \frac{a_4}{2} 
r_2 \cos \chi - \varepsilon^2\left(\frac{1}{24}a_1^2r_1^2 +\frac{1}{2}a_1a_2r_2^2+ e^{-2 \tau}(\frac{1}{8} a_3a_4r_2^2 
+\frac{1}{64} a_4^2(9r_1^2+4r_2^2)\right) ,\\

\dot{\psi}_2 & =  - \varepsilon e^{- \tau} \frac{a_4}{8} \frac{r_1^2}{r_2} \cos \chi -\varepsilon^2\left(\frac{1}{4}a_1a_2r_1^2 
+\frac{1}{30}a_2^2r_1^2+ \frac{29}{120}a_2^2r_2^2 +e^{-2 \tau}(\frac{1}{16}a_3a_4r_1^2+\frac{1}{32}a_4^2r_1^2 
+\frac{5}{96}a_3^2r_2^2)\right). 
\end{cases} 
\end{eqnarray} 

The result is surprising as we would expect $\tau$-dependent terms for the amplitudes at second order; such terms 
arise only for the angles. Also combination angles for $\psi_1, \psi_2$ are not present at this level of averaging-normalisation. 
The 2nd order system \eqref{ave2} was computed without the time-independence in \cite{TV00}, eq. (4.2). 
Leaving out the time-dependent terms the results agree. \\
The system has interesting implications: 

\begin{enumerate}

\item During an interval of time of order $1/ \varepsilon$ the integrals (adiabatic invariants) \eqref{intave1} are active,  
will govern the orbital dynamics and accordingly the corresponding distribution function in phase-space. This holds for 
$n= 1, 2, \ldots$

\item If $\delta= \varepsilon$ and on asymptotically longer time intervals like $1/ \varepsilon^2$ the time-independent 
system involving the coefficients $a_1, a_2$ dominates the dynamics. In \cite{TV00} it is shown that for 
this system, depending on $a_1, a_2$, 2 resonance manifolds can exist on the energy manifold. Introducing 
the combination angle: 
\begin{equation} \label{resangle} 
\chi_2 = 4 \psi_1 - 2 \psi_2, 
\end{equation} 
we find according to \cite{TV00}  that on intervals of time larger than $1/ \varepsilon$ we have: 
\begin{equation} \label{resmaneq}
\frac{d \chi_2}{dt} = \varepsilon^2 \left((- \frac{1}{6}a_1^2 + \frac{1}{2}a_1a_2 + \frac{1}{15}a_2^2)r_1^2  + (- 2a_1a_2 + 
\frac{29}{60}a_2^2) r_2^2\right).
\end{equation} 
Resonance manifolds exist if the righthand side of eq. \eqref{resmaneq} has a zero, the combination angle $\chi_2$ is 
not timelike. In this case the resonance manifold 
with $4 \psi_1 - 2 \psi_2 =0$ has stable $2:4$ resonant periodic orbits surrounded by tori, 
for $4 \psi_1 - 2 \psi_2 = \pi$ the $2:4$ resonant 
periodic orbits also exist but are unstable. The resonance manifolds have size $O(\varepsilon)$, the dynamics takes 
place on intervals of time of order $1/ \varepsilon^3$; for details see \cite{TV00}.  \\
Outside the resonance manifolds the dynamics is characterised by the quadratic integrals $E_1, E_2$ for each of 
the 2 dof.  
For the interaction of the 2 modes on these long time intervals coefficient $a_2$ is essential, but 
note that if $a_1=0$, the resonance manifolds do not exist as the righthand side of eq..~\eqref{resmaneq} has no zero. 
\end{enumerate}

\subsection{Consequences for the distribution function}  
Suppose we start with a collection of particles (stars) characterised by a distribution function satisfying the 
Boltzmann equation that is nearly collisionless as it has small dynamical friction added. The system is already in 
axi-symmetric state but the evolution to mirror symmetry to the galactic plane is still going on. 
On an interval of time $O(1/ \varepsilon)$, the first stage of evolution, we have, apart from $J$, 2 active integrals of motion: 
$ E_0, I_3$ and a resonance manifold for each value of the energy described by eq. \eqref{resman1}. The family of periodic  
solutions with constant amplitude (constant to $O(\varepsilon^2$)) on the energy manifold will be stable for 
combination angle $\chi = 0$.
The distribution function will be a function of the 3 integrals. The velocity distribution $v_1, v_2$ and their dispersion will 
depend on the existence of these integrals. 

\begin{figure}
\begin{center}
\resizebox{16cm}{!}{
\includegraphics{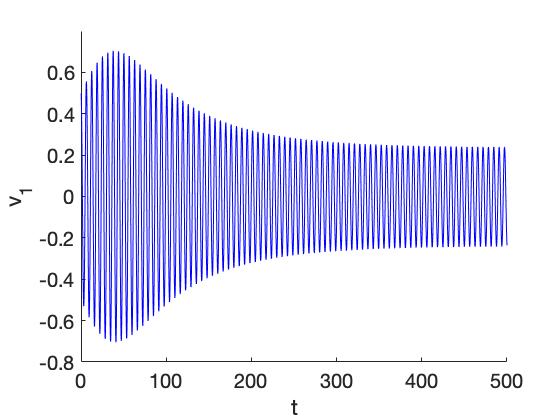} \, \includegraphics{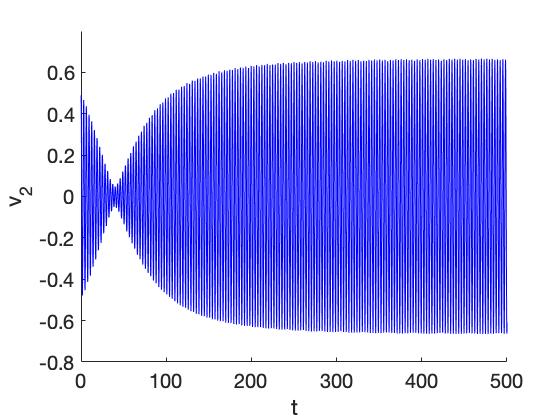}} 
\resizebox{16cm}{!}{
\includegraphics[width=8cm]{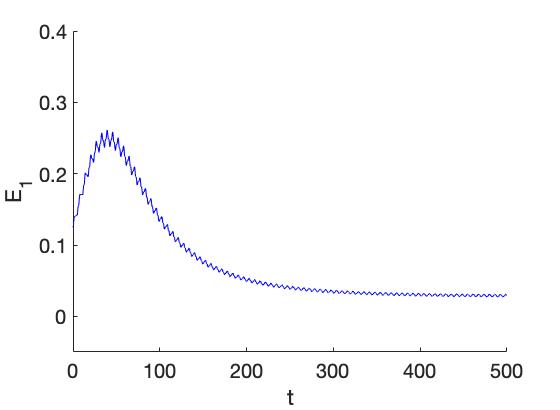} \, \includegraphics[width=8cm]{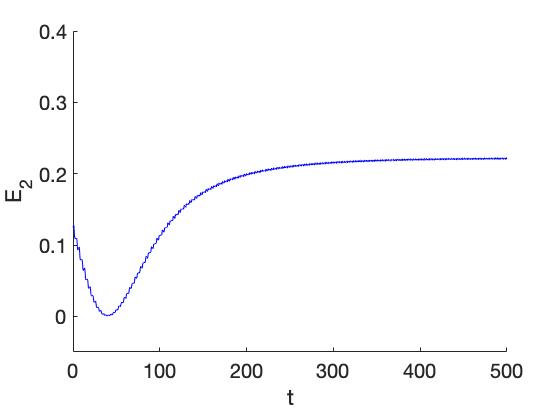}}
\caption{The behaviour of $v_1(t), v_2(t)$ of system \eqref{eqs} in the case $n=2$. We have $\alpha(\delta t)= \exp{(- \varepsilon^n t)}, 
\varepsilon =0.1, \omega =2, a_1= a_2= 1, a_3= 0.75, a_4=1.5$ with initial conditions $q_1(0)=  q_2(0)=0, v_1(0)=v_2(0)
=0.5$. After around 100 timesteps the time-dependent interaction vanishes, after 200 timesteps the velocities have 
become different. Below the actions with $E_1(0)= E_2(0)= 0.125$.
\label{fig1}}
\end{center}
\end{figure}

On intervals of time asymptotically larger than $O(1/ \varepsilon)$, the primary 
resonance manifolds described by eq. \eqref{resman1} vanish, they are replaced by the smaller resonance manifolds 
(of size $\varepsilon$) located by 
the zeros of eq. \eqref{resmaneq}  and $\chi_2 = 0, \pi$. The distribution function will now evolve to a function of 
$ E_1, E_2$. See fig.~ \ref{fig1} where $n=2$ so that the symmetric state develops on intervals of time $O(1/ \varepsilon^2)$.  
On an interval of time order $1/ \varepsilon^2$ there is first 
a typical $1:2$ resonance exchange between the 2 dof after which the dynamics settles at a slightly lower amplitude. 

\begin{figure}
\begin{center}
\resizebox{16cm}{!}{
\includegraphics[width=8cm]{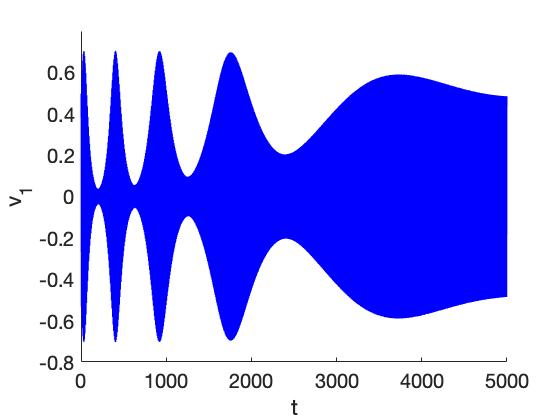} \, \includegraphics[width=8cm]{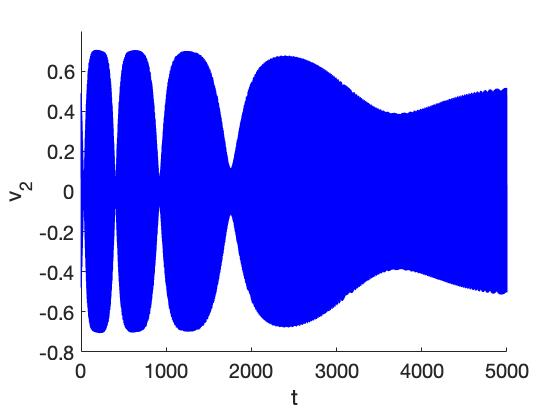}}\\
\resizebox{16cm}{!}{
\includegraphics[width=8cm]{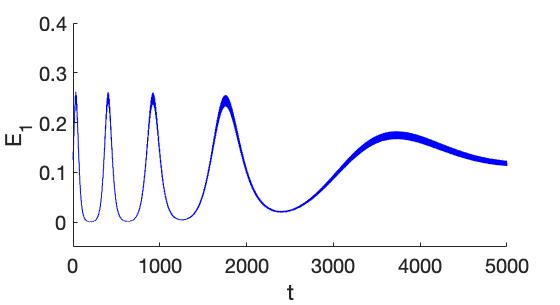} \, \includegraphics[width=8cm]{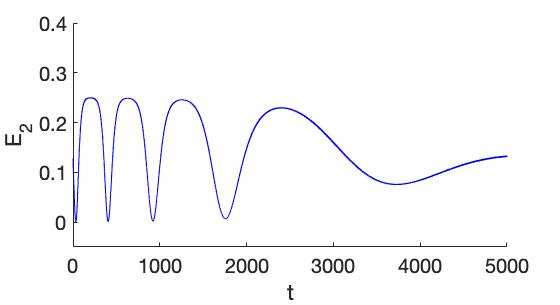}}
\caption{The behaviour of $v_1(t), v_2(t)$ of system \eqref{eqs} in the case $n=3$, it takes longer for the time-dependent 
interaction to vanish. Except for the value of $n$ the parameter values and initial 
conditions are as in fig.~\ref{fig1}. If $n=3$ it takes at least 5000 timesteps to stabilise the dynamics. 
Below are the corresponding actions. \label{fig2}}
\end{center}
\end{figure}

In fig.~\ref{fig2} we have $n=3$ so that on intervals of size $O(1/ \varepsilon^2)$ we have still resonant interaction, the 
symmetric state develops on intervals of size $O(1/ \varepsilon^3)$. 
The exchanges and time evolution are first much more frequent as it takes longer for the 
asymmetric terms to vanish. Then the system settles at a much 
lower oscillation amplitude, in a sense experiencing more of its asymmetric past. \\

When the system is close to mirror symmetry, we will find for each value of the energy  in a resonance manifold 
 a family of tori around the stable periodic solutions, 
so for varying energy values this will be a 2-parameter family. Outside these resonance manifolds the orbits will move 
quasi-linearly, only the phases are position dependent. 

The choice of $\delta$, the parameter 
determining the timescale of destroying the asymmetry of the force field, together with the energy level 
will determine the resulting positions and velocities. 

\section{The $1:3$ resonance} 
The $1:3$ resonance is for Hamiltonian \eqref{Ham} dynamically very different. Most of the analysis can be 
deduced from  \cite{FV79}, 
we will summarise the results. Putting $\omega =3$ in \eqref{Ham}  we find after 2nd order averaging: 
\begin{equation} \label{ampl13} 
\dot{r}_1 = O(\varepsilon^3),\, \dot{r}_2 = O(\varepsilon^3). 
\end{equation} 
The result is quite remarkable as for the $1:3$ resonance the slowly vanishing asymmetry of the potential plays no 
part for the amplitudes  to order 3 in $\varepsilon$. For the phases we find the same phenomenon but with nontrivial variation: 
\begin{eqnarray} \label{phase13} 
\dot{\psi}_1 & = - \varepsilon^2 \left(\frac{5}{12}a_1^2r_1^2 + (\frac{1}{2}a_1a_2 - \frac{1}{35}a_2^2)r_2^2\right),\\
\dot{\psi}_2 & = - \varepsilon^2 \left( (\frac{1}{6}a_1a_2+ \frac{1}{105}a_2^2)r_1^2 + \frac{23}{140}a_2^2r_2^2\right).
\end{eqnarray} 
if $\delta =O(\varepsilon)$ or $O(\varepsilon^2)$ the asymmetry plays no significant part. 
The theory of higher order resonance, see \cite{SVM} and for extension and examples \cite{TV00}, shows that the 
combination angle $\chi_3= 6 \psi_1 -2 \psi_2$ plays a crucial role. The equation for $\chi_3$ becomes: 
\begin{equation} \label{resangle13} 
\frac{d \chi_3}{dt} = -  \varepsilon^2 \left( (\frac{5}{2}a_1^2 - \frac{1}{6}a_1a_2- \frac{1}{105}a_2^2)r_1^2 - 
(3 a_1a_2 + \frac{47}{140})r_2^2 \right). 
\end{equation} 
Zero solutions of the righthand side of eq.~\eqref{resangle13}  produce, together with the condition $\chi_3= 0, \pi$ 
resonance manifolds of size $O(\varepsilon^2)$ with characteristic timescale $O(1/ \varepsilon^4)$ for the orbits and 
motion on the corresponding tori; these estimates follow from the analysis in \cite{TV00}. \\
It is clear that if $a_1=0$, no resonance manifolds of this type are present. 

\section{The $1:1$ resonance} 
If the epicyclic frequency and the vertical frequency are equal or close, the $1:1$ resonance becomes important. With  
$\alpha(\delta t)$ given by eq.~\eqref{alfa} the equations of motion induced by Hamiltonian \eqref{Ham} iare:
\begin{eqnarray} \label{eqs11}
\begin{cases}
\ddot{q}_1 + q_1 & =   \varepsilon (a_1q_1^2 +a_2q_2^2) + \varepsilon  \alpha(\delta t) 2a_4q_1q_2,\\ 
\ddot{q}_2 +  q_2 & =   \varepsilon 2a_2q_1q_2 + \varepsilon  \alpha(\delta t) (a_3q_2^2 +a_4q_1^2).
\end{cases} 
\end{eqnarray}  
By averaging system \eqref{slow} if $\omega=1$ 
we find that all first order averaged terms vanish. Significant dynamics takes place on a longer timescale; we choose 
$n \geq 2$  when considering longer time scales. 
Second order averaging based on \cite{SVM} produces with 
$\chi = \psi_1 - \psi_2$ the system: 
\begin{eqnarray} \label{ave3}
\begin{cases} 
\dot{r}_1 & =  \varepsilon^2 (\frac{1}{12}a_1a_2 - \frac{1}{2}a_2^2)r_1r_2^2 \sin 2 \chi + \varepsilon^2 \alpha(\delta t) 
(\frac{1}{12}a_3a_4 - \frac{1}{2}a_4^2)r_1r_2^2 \sin 2 \chi , \\
\dot{r}_2 & = - \varepsilon^2 (\frac{1}{12}a_1a_2 - \frac{1}{2}a_2^2)r_1^2r_2 \sin 2 \chi - \varepsilon^2 \alpha(\delta t) 
(\frac{1}{12}a_3a_4 - \frac{1}{2}a_4^2)r_1^2r_2 \sin 2 \chi , \\
\dot{\psi}_1 & = - \varepsilon^2\left( \frac{5}{12}a_1^2r_1^2+ (\frac{1}{2}a_1a_2+ \frac{1}{3}a_2^2)r_2^2 -(\frac{1}{12}a_1a_2 
- \frac{1}{2}a_2^2)r_2^2 \cos 2 \chi \right) - \\ 
& \hspace{0.5 cm} \varepsilon^2  \alpha(\delta t) \left((\frac{1}{2}a_3a_4 + \frac{1}{3}a_4^2)r_2^2 + \frac{5}{12}a_4^2r_1^2 -
(\frac{1}{12}a_3a_4- \frac{1}{2}a_4^2)r_2^2 \cos 2 \chi \right), \\
\dot{\psi}_2 & = - \varepsilon^2\left( (\frac{1}{2}a_1a_2+ \frac{1}{3}a_2^2)r_1^2 + \frac{5}{12}a_2^2r_2^2 - (\frac{1}{12}a_1a_2 
- \frac{1}{2}a_2^2)r_1^2 \cos 2 \chi \right) - \\ 
& \hspace{0.5 cm} \varepsilon^2  \alpha(\delta t) \left( \frac{5}{12}a_3^2r_2^2 +(\frac{1}{2}a_3a_4+ \frac{1}{3}a_4^2)r_1^2 
-(\frac{1}{12}a_3a_4 - \frac{1}{2}a_4^2)r_1^2 \cos 2 \chi \right). 
\end{cases}
\end{eqnarray} 
From system \eqref{ave3} we can derive the equation for $\chi$. The solutions of the system with appropriate initial values 
produce an $O(\varepsilon^2)$ approximation on an interval of time order $1/ \varepsilon$ of system \eqref{slow} with 
$\omega =1$ but an $O(\varepsilon)$ approximation on  an interval of time order $1/ \varepsilon^2$. The second estimate 
is a kind of ``trade-off'' of error estimates valid under special conditions formulated in \cite{B76}; see also \cite{SVM}. 
It is remarkable that system \eqref{ave3} shows full resonance involving exchange of energy between the 2 dof even when 
the asymmetric potential terms have become negligible.  \\
A systematic study of the symmetric case ($\alpha(\delta t)=0, t \geq 0$) can be found in \cite{FV79}. We summarise the 
results for the symmetric case: 
\begin{itemize} 
\item System \eqref{ave3} has in the symmetric case 2 integrals of motion: 
\begin{equation} \label{11Int1}
E_0= \frac{1}{2}(r_1^2 +r_2^2) = \frac{1}{2}(\dot{q}_1^2 +q_1^2 + \dot{q}_2^2 +q_2^2)
\end{equation} 
and (for historical reasons called $I_3$): 
\begin{equation} \label{11Int2} 
I_3= r_1^2r_2^2 \cos 2 \chi + \alpha r_1^4 + \beta r_1^2, 
\end{equation}
with $\alpha, \beta$ rational functions of $a_1, a_2$. We leave out some degenerate cases of the coefficients, see 
\cite{FV79}. 

\item The $q_1, \dot{q}_1$ normal mode is an exact solution, it is unstable for $- 1/3 < a_1/(3a_2) < 2/15$ and 
$1/3<  a_1/(3a_2) < 2/3$.
\item The $q_2, \dot{q}_2$ normal mode is obtained from the system \eqref{ave3} as an $O(\varepsilon)$ approximation. 
It is unstable for   $- 1/3 < a_1/(3a_2) < 1/3$. 
\item The in-phase periodic solutions $\chi = 0, \pi$ exist for $ a_1/(3a_2) < 2/3$ and are stable for $-1/3 <  a_1/(3a_2) < 2/3$. 
\item The out-of-phase periodic solution $\chi = \pi/2, 3 \pi/2$ exist and are stable for $ a_1/(3a_2) < 2/15$. 
\end{itemize} 
We expect the dynamics of the symmetric case to describe the orbits of system \eqref{eqs11} on intervals of time larger 
than $1/ \varepsilon^n$. The dynamics will then be governed by the integrals \eqref{11Int1} and \eqref{11Int2} producing 
a complicated velocity distribution and varying actions. \\
On the long (starting) interval of order $1/ \varepsilon^2$ system \eqref{ave3} has also 2 integrals. 
Remarkably enough integral \eqref{11Int1} holds for all time $t \leq 0$. 



\section{Appendix} \label{app} 
We present a modification of the averaging technique.  
Consider the $T$-periodic vector fields $f_1, f_2$ and the slowly varying ODE: 
\begin{equation} \label{slow2} 
\dot{x} = \varepsilon f_1(t, x) + \varepsilon f_2(t, x),
\end{equation} 
With $x, f_1 , f_2 \in   {\mathbb{R}}^n$ and $f_1, f_2$ twice continously differentiable in a bounded domain $D$ in 
${\mathbb{R}}^n$, continuously differentiable in $t$. Suppose in addition that:
\[ \frac{1}{T} \int_0^T f_1(t, x)dt =0, \]
where $x$ is kept constant during integration. We use the near-identity transformation:
\begin{equation} \label{nearid} 
x(t) = y(t) + \varepsilon u(t, y(t)),\, u(t, y(t)) = \int_0^t f_1(s, y(t)) ds.
\end{equation} 
As the vector field $f_1$ is $t$-periodic with average zero, $u(t, y(t))$ is bounded in $D$ by a constant independent of 
$\varepsilon$. Substituting $x(t)$ with eq. \eqref{nearid} in ODE \eqref{slow2} we have: 
\[ \dot{x}= \dot{y} + \varepsilon f_1(t, y) + \varepsilon \frac{\partial u}{\partial y} \dot{y} = \varepsilon f_1(t, y+ \varepsilon u) + 
\varepsilon f_2(t, y+ \varepsilon u).  \]
We expand $\varepsilon f_1(t, y+ \varepsilon u) + \varepsilon f_2(t, y+ \varepsilon u)= \varepsilon f_1(t, y) + 
\varepsilon f_2(t, y) + O(\varepsilon^2)$. and find: 
\[ \left(I  + \varepsilon \frac{\partial u}{\partial y}\right) \dot{y} = \varepsilon f_2(t, y) + O(\varepsilon^2). \] 
$I$ is the $n \times n$ unit matrix, the matrix $I  + \varepsilon \partial u / \partial y$ has a bounded inverse, so that: 
\begin{equation} \label{transf}
 \dot{y} = \varepsilon f_2(t, y) + O(\varepsilon^2). 
 \end{equation} 
 where the $O(\varepsilon^2)$ terms can be computed explicitly. The procedure is removes the non-resonant 
 terms from the righthand side of eq. \eqref{slow2}; this is useful if we are able to perform analysis on the resonant 
 part with explicit slow time as in section \ref{sec3}. 

\medskip 
{\bf Acknowledgement}\\
System \eqref{ave2} was computed by Taoufik Bakri using {\sc Mathematica}.

\end{document}